\documentclass[11pt,letterpaper]{amsart}
\usepackage{amsmath,amssymb}
\usepackage{mathrsfs}
\usepackage{amsfonts}
\usepackage[dvipsnames]{xcolor}
\usepackage{graphicx}
\usepackage{float}  
\usepackage{cite}
\usepackage{microtype}
\usepackage{cases}
\usepackage{comment}
\usepackage{indentfirst}
\allowdisplaybreaks
\theoremstyle{plain}
\newtheorem{theorem}{Theorem}[section]
\newtheorem{lemma}[theorem]{Lemma}
\newtheorem{proposition}[theorem]{Proposition}

\theoremstyle{definition}

\theoremstyle{remark}

\numberwithin{equation}{section}

\DeclareMathOperator{\tr}{tr}
\DeclareMathOperator{\dist}{dist}
\newcommand{\R}{\mathbb R}
\newcommand{\Om}{\Omega}
\newcommand{\eps}{\varepsilon}
\newcommand{\Ga}{\Gamma_2}
\newcommand{\cD}{\mathcal D}
\newcommand{\cT}{\mathcal T}
\newcommand{\cK}{\mathcal K}
\newcommand{\cF}{\mathscr F}
\newcommand{\cL}{\mathcal L}
\newcommand{\dm}{\,d\mu}
\newcommand{\norm}[1]{\left\|#1\right\|}
\newcommand{\abs}[1]{\left|#1\right|}

\usepackage[colorlinks=true,citecolor=blue,linkcolor=blue]{hyperref}

\begin{document}
\title[Global $W^{2,p}$ and $C^{2,\alpha}$ regularity]
{Global $W^{2,p}$ and $C^{2,\alpha}$ regularity for the sigma-$2$ equation}

\author[Y. Lyu]{Yasheng Lyu}
\address{School of Mathematics and Statistics, Xi'an Jiaotong University, Xi'an, Shaanxi 710049, People's Republic of China}
\email{lvysh21@stu.xjtu.edu.cn}

\date{}
\keywords{$2$-Hessian equation, global $W^{2,p}$ estimates, global $C^{2,\alpha}$ estimates, Hessian measures, viscosity solutions.}

\begin{abstract}
In this paper, we prove that continuous $2$-convex solutions of the $\sigma_{2}$ equation with positive continuous density on bounded strictly mean-convex $C^{3}$ domains with $C^{3}$ boundary values belong to $W^{2,p}$ for every $1<p<\infty$ in every dimension.
If, in addition, the density is $C^{\alpha}$, $0<\alpha<1$, the solutions belong to $C^{2,\alpha}$.
Moreover, two counterexamples show that neither the $C^{3}$ assumption on the domain nor that on the boundary values can in general be weakened to $C^{2,1}$.
\end{abstract}

\maketitle

\section{Introduction}\label{sec:intro}

We consider the Dirichlet problem
\begin{equation}\label{eq:dirichlet}
\sigma_2(D^2u)=f\quad\text{in }\Om,\qquad
u=\varphi\quad\text{on }\partial\Om,
\end{equation}
where $\Om\subset\R^n$, $n\ge3$, is a bounded domain. For symmetric matrices,
we write
\[
\sigma_2(A)=\frac12\bigl[(\tr A)^2-|A|^2\bigr],\qquad
\Ga=\{A:\tr A>0,\ \sigma_2(A)>0\}.
\]
Here $|A|$ is the Frobenius norm and $\Ga$ is the ellipticity cone.
A continuous function is $2$-convex
if every $C^2$ function touching it from above has Hessian in
$\overline\Ga$ at the contact point; this class is denoted by $\Phi_2(\Om)$.
We interpret \eqref{eq:dirichlet} through the Hessian measures of
Trudinger and Wang~\cite{TW1,TW2}:
\[
\mu_2[u]=f\,dx,
\]
where $\mu_2[u]=\sigma_2(D^2u)\,dx$ for smooth $u$.
For continuous positive $f$, this agrees with the admissible viscosity
formulation; see \cite[Section~2]{TW2}.

We use the curvature convention
$\mathrm{II}(\tau,\tau)=D_\tau\nu\cdot\tau$ for the outward normal $\nu$.
Thus spheres have positive principal curvatures, and we write
$H_{\partial\Om}=\tr\mathrm{II}$.
The condition $H_{\partial\Om}>0$ is the classical domain condition
for the $2$-Hessian Dirichlet problem of Caffarelli, Nirenberg and
Spruck~\cite{CNS}. Trudinger~\cite{Trudinger} extended boundary
second-derivative estimates to Hessian quotient equations.
Ivochkina, Trudinger and Wang~\cite{ITW} developed boundary
second-derivative estimates for degenerate Hessian equations using
vector fields with skew-symmetric Jacobians.

\begin{samepage}
\begin{theorem}\label{thm:wp}
Let $\Om\subset\R^n$, $n\ge3$, be a bounded strictly mean-convex domain
with $\partial\Om\in C^3$. Suppose $\varphi\in C^3(\partial\Om)$ and
$f\in C(\overline\Om)$, with $0<\lambda\le f\le\Lambda$.
Then any $u\in C(\overline\Om)\cap\Phi_2(\Om)$ solving
\eqref{eq:dirichlet} in the Hessian-measure sense belongs to
$W^{2,p}(\Om)$ for every $1<p<\infty$ and satisfies
\begin{equation}\label{eq:main-estimate}
\|u\|_{W^{2,p}(\Om)}\le C_p,
\end{equation}
where $C_p$ depends only on $n,p,\lambda,\Lambda,\Om,\varphi$
and the modulus of continuity of $f$.
\end{theorem}
\end{samepage}

\begin{theorem}\label{thm:c2alpha}
Let $\Om\subset\R^n$, $n\ge3$, be a bounded strictly mean-convex domain
with $\partial\Om\in C^3$. Suppose $\varphi\in C^3(\partial\Om)$ and
$f\in C^\alpha(\overline\Om)$ for some $0<\alpha<1$, with $f\ge\lambda>0$.
Then any $u\in C(\overline\Om)\cap\Phi_2(\Om)$ solving
\eqref{eq:dirichlet} in the Hessian-measure sense belongs to
$C^{2,\alpha}(\overline\Om)$ and satisfies
\begin{equation}\label{c2:eq:main-estimate}
\|u\|_{C^{2,\alpha}(\overline\Om)}\le C,
\end{equation}
where $C$ depends only on $n,\alpha,\lambda$,
$\|f\|_{C^\alpha(\overline\Om)}$, $\Om$ and $\varphi$.
\end{theorem}

The dependence on $\Om$ and $\varphi$ includes moduli of continuity of
the third derivatives of boundary charts and boundary values. This dependence
is uniform on compact families in the $C^3$ topology with a common
positive lower bound for the mean curvature.
The $C^3$ assumptions cannot in general be replaced
by $C^{2,1}$. Section~\ref{sec:examples} gives density-one solutions
with either a smooth domain and $C^{2,1}$ boundary values, or a uniformly
convex $C^{2,1}$ domain and polynomial boundary values. In both examples,
$W^{2,p}$ regularity holds exactly for $p<n+1$.

For constant density, Warren and Yuan~\cite{WY} proved interior Hessian
estimates in dimension three, and Shankar and Yuan~\cite{SY4} obtained
the unrestricted estimate in dimension four. For variable right-hand
sides, Qiu~\cite{Qiu} proved the three-dimensional estimate for positive
$C^2$ data. In dimension three, Zhou~\cite{ZhouGSL} obtained interior
Hessian estimates and $C^{2,\alpha}$ regularity for continuous
admissible viscosity solutions with positive Lipschitz density, while
Xu~\cite{Xu} obtained perturbative interior $C^{2,\alpha}$ estimates
under a smallness assumption on the H\"older seminorm of the density.
Fan~\cite{Fan} treated positive $C^{1,1}$ data in dimension four.
In arbitrary dimension, Guan and Qiu~\cite{GQ}
obtained interior estimates under a lower bound for $\sigma_3(D^2u)$,
in particular for convex solutions. McGonagle, Song and Yuan~\cite{MSY}
obtained the constant-density estimate under an almost-convexity
condition, and Shankar and Yuan~\cite{SYsemi} treated smooth semiconvex
solutions. Interior regularity for convex and almost-convex viscosity
solutions of the constant equation was proved by Mooney~\cite{MooneyConvex}
and Shankar and Yuan~\cite{SYalmost}, respectively.

For positive Lipschitz density, Chen, Jian and Zhou~\cite{CJZ}
obtained interior Hessian estimates for smooth convex solutions, and
Chen, Jian, Tu and Zhou~\cite{CJTZ} proved interior $C^2$ regularity
for convex viscosity solutions.
Zhou and Zhu~\cite{ZZ} obtained interior $C^{2,\alpha}$ regularity for
convex viscosity solutions with positive $C^\alpha$ density.
Li and Wu~\cite{LW} established interior
Hessian estimates for the constant equation on the full positive branch
in every dimension; their revised version also treats positive
Lipschitz right-hand sides $f(x,u)$. Chen, Zhou and Zhu~\cite{CZZ}
proved interior $C^{2,\alpha}$ estimates and regularity on that branch
for positive $C^\alpha$ density.

Interior Hessian estimates from integral bounds were obtained by
Urbas~\cite{Urbas,Urbas2001}; Dinew and Ko{\l}odziej~\cite{DK}
improved the integrability threshold for a range of real Hessian
equations. For the constant $\sigma_2$ equation,
Mooney~\cite{Mooney} derived an interior Hessian estimate from
$W^{2,p}$ for every $p>2$. For densities bounded above and below by
positive constants, the author~\cite{Lyu} proved interior $W^{2,1}$
regularity of continuous $2$-convex solutions; if the density is
continuous, the same work gives interior $W^{2,p}$ estimates for every
finite $p>1$.

For the Monge--Amp\`ere equation, Caffarelli~\cite{Caffarelli} proved
interior $W^{2,p}$ estimates for strictly convex solutions with positive
continuous density. Savin~\cite[Theorem~1.2 and Corollary~1.4]{Savin}
obtained global $W^{2,p}$ estimates under quadratic separation at the
boundary, in particular on uniformly convex $C^3$ domains with $C^3$
Dirichlet data. Trudinger and Wang~\cite[Theorem~1.1]{TW} proved global
$C^{2,\alpha}$ estimates under the latter boundary assumptions for
positive $C^\alpha$ density. Savin~\cite[Theorem~1.1]{SavinC2}
subsequently obtained pointwise boundary $C^{2,\alpha}$ estimates under
quadratic separation, with pointwise $C^\alpha$ density and pointwise
$C^{2,\alpha}$ boundary and boundary values. Wang's examples~\cite{Wang}
show that $C^{2,1}$ regularity of either the domain or the boundary
data does not suffice for all finite $W^{2,p}$ exponents; see also
\cite{Savin,TW}. In dimension two, where $\sigma_2=\det$, the
corresponding statements of Theorems~\ref{thm:wp}
and~\ref{thm:c2alpha} are contained in these results.

The boundary part of the $W^{2,p}$ proof couples a weighted quadratic
comparison with simultaneous control of the trace and normalized
tangential gradients. It uses the comparison framework of~\cite{CZZ}
and the trace-tail argument of~\cite{Lyu}, with a positive lower bound
for the tangential trace established under the present $C^3$
assumptions. For $C^\alpha$ density,
the global $W^{2,p}$ bound and the second-difference and resolvent
estimates of~\cite{CZZ} are supplemented by additive comparisons that
give uniform bounds for the boundary mixed derivatives and the
global trace.

\section{Global \texorpdfstring{$W^{2,p}$}{W2,p} estimates}\label{sec:wp}

For symmetric matrices write $A:B=\tr(AB)$, and set
\[
G(A)=\sqrt{\sigma_2(A)},\qquad N(A)=(\tr A)I-A.
\]
On $\Ga$, $N(A)>0$ and $|A|<\tr A$; $G$ is concave and homogeneous
of degree one. Fix $0<m<M$ and write
$h_0=\min_{\partial\Om}H_{\partial\Om}>0$.
Throughout this section, $C,c$ depend only on $n,m,M$ and the geometric
and boundary data in Theorem~\ref{thm:wp}. The constant $C_b$ may additionally depend
on $\|b\|_{C^2}$, but never on the truncation parameters or derivatives of $g$.

\subsection{A lower bound for the tangential trace}
\label{wp:nondegeneracy}

There is a defining function $\rho<0$ such that, with $r=\dist(x,\partial\Om)$,
\begin{equation}\label{eq:rho}
\rho=0\ \text{on }\partial\Om,\qquad -\rho\asymp r,\qquad
D^2\rho-\theta I\in\Ga\quad\text{on }\overline\Om
\end{equation}
for a fixed $\theta>0$. We recall the construction to specify its regularity.
In a collar take $\rho_0=-r+Jr^2$. Its Hessian at the boundary has eigenvalues
$(\kappa_1,\ldots,\kappa_{n-1},2J)$, so
\[
\sigma_2(D^2\rho_0)=\sigma_2(\kappa)+2JH_{\partial\Om}>0
\]
for large fixed $J$. Shrink the collar to $0\le r\le\ell$, with
$J\ell\le1/8$. Glue $\rho_0$ to a small strictly convex quadratic whose
range lies in $(-\ell/2,-\ell/4)$, using a smooth maximum in the collar
with smoothing scale less than $\ell/32$.
This equals $\rho_0$ near the physical boundary and the quadratic near the
inner edge. Its Hessian is a convex combination of the two Hessians plus a
nonnegative matrix. Compactness gives \eqref{eq:rho}; the construction is
$C^3$ for $C^3$ boundaries and smooth for smooth boundaries.

Set $d=-\rho$. For any admissible $v$, let
$a=DG(D^2v)$ and $\cL=a^{ij}\partial_{ij}$. Concavity and homogeneity give
\[
a:D^2\rho\ge G(D^2\rho-\theta I)+\theta\tr a\ge\theta\tr a.
\]
Consequently the fixed defining function satisfies
\begin{equation}\label{eq:linear-barrier}
\cL d\le-\theta\tr a,\qquad
N(D^2\rho)\ge(n-1)\theta I.
\end{equation}

Take a fixed ambient $C^3$ extension $\Phi$ of $\varphi$, smooth when
the domain and boundary data are smooth. Such an extension is obtained
in boundary graph coordinates and patched by a smooth partition of unity.
Choose $A_0$ large enough that
\[
\begin{gathered}
D^2(\Phi+A_0\rho)\in\Ga,\qquad
G(D^2(\Phi+A_0\rho))\ge M,\\
\Delta(\Phi-A_0\rho)\le0.
\end{gathered}
\]
These conditions follow from \eqref{eq:rho}. If $v=\varphi$ on
$\partial\Om$ and $m\le G(D^2v)\le M$, comparison with $\Phi+A_0\rho$ and the maximum principle for
$v-(\Phi-A_0\rho)$ give
\begin{equation}\label{eq:comparison}
\begin{gathered}
\Phi+A_0\rho\le v\le\Phi-A_0\rho,\qquad |v-\Phi|\le A_0d,\\
\|v\|_{L^\infty(\Om)}+\|Dv\|_{L^\infty(\partial\Om)}\le C.
\end{gathered}
\end{equation}
The boundary gradient bound follows by taking normal derivatives of the
comparison inequalities.

The following boundary estimate uses only the positive lower bound
for the density.

\begin{proposition}\label{prop:trace}
Let $\Om\subset\R^n$, $n\ge3$, be a bounded $C^3$ domain with
$H_{\partial\Om}\ge h_0>0$, and let $\varphi\in C^3(\partial\Om)$.
If
$v\in C^2(\overline\Om)$ satisfies $D^2v\in\Ga$ in $\Om$,
$v=\varphi$ on $\partial\Om$,
and $\sigma_2(D^2v)\ge m^2$ in $\Om$ for some $m>0$, then
\begin{equation}\label{eq:tracegap}
\tr\bigl(D^2v|_{T\partial\Om}\bigr)\ge c_0>0.
\end{equation}
Here $c_0$ depends only on $n,m,h_0$, the quantitative $C^3$ boundary
geometry and data, and their third-derivative moduli of continuity.
\end{proposition}

\begin{proof}
Fix a boundary point and choose coordinates so that it is the origin,
the inward normal is $e_n$, and locally $\Om=\{x_n>\beta(x')\}$, with
$\beta(0)=D\beta(0)=0$. Subtract a tangential affine function
so that the boundary value $\phi(x')$ has vanishing constant and linear
terms. Write
\[
\beta=\beta_2+\beta_3+o(s^3),\qquad
\phi=\phi_2+\phi_3+o(s^3),\qquad s=|x'|,
\]
where subscripts denote homogeneous Taylor terms. Put
\[
H_0=\Delta'\beta_2\ge h_0,\qquad q_0=\frac{\Delta'\phi_2}{H_0},
\qquad Q_2=\phi_2-q_0\beta_2,\quad P_3=\phi_3-q_0\beta_3.
\]
Thus $\Delta'Q_2=0$. Choose $0<\tau<h_0/[4(n-1)]$, and set
\[
d_0=x_n-\beta_2(x')+\frac{\tau}{2}s^2,\qquad
e=\frac{Dd_0}{|Dd_0|},\qquad
\cT=(I-e\otimes e):D^2.
\]
In the small cap
\[
U_h=\{s<h,\ \beta(x')<x_n<h^2\}
\]
we have
\begin{equation}\label{eq:cap}
d_0\asymp s^2+|x_n|,\qquad \cT d_0\le-c<0,
\qquad d_0=\tau s^2/2+O(s^3)\quad\text{on }\partial\Om.
\end{equation}
Indeed, $d_0\ge x_n-\beta(x')+\tau s^2/4$ for small $h$, and
$\cT d_0=-H_0+(n-1)\tau+O(s^2)$.

Define the cubic polynomial
\[
\ell_1=\frac{\Delta'P_3}{H_0+2\tau},\qquad
w_0=Q_2+q_0x_n+P_3+\ell_1(x')(d_0-\tau s^2/2).
\]
It satisfies
\begin{equation}\label{eq:model}
\begin{gathered}
|\cT w_0|\le Cd_0,\qquad |D^2w_0|\le C,\\
|w_0-\phi|\le\epsilon(h)s^3\quad\text{on }\partial\Om,
\qquad \epsilon(h)\longrightarrow0.
\end{gathered}
\end{equation}
To check the tangential trace, observe that
$(I-e\otimes e)Dd_0=0$ and $D'd_0=O(s)$. Thus
\[
\cT(\ell_1d_0)=\ell_1\cT d_0,\qquad
\cT(\ell_1s^2/2)=(n+1)\ell_1+O(s^3),
\]
while $\cT Q_2=O(s^2)$ and
$\cT P_3=\Delta'P_3+O(s^3)$. Consequently,
\[
\cT w_0=\Delta'P_3-(H_0+2\tau)\ell_1+O(s^2)=O(d_0).
\]
The boundary error is $o(s^3)+O(s^4)$, since
$d_0-\tau s^2/2=\beta-\beta_2$ on $\partial\Om$.

Choose fixed $C_0$ with $C_0s^2-\beta(x')\ge s^2$, and set
\[
V=w_0+\frac{B_0}{2}d_0^2+\eta_0(C_0s^2-x_n).
\]
For $B_0\ge B_*$, with $B_*\ge1$ fixed, and $0<\eta_0\le1$,
\eqref{eq:cap}--\eqref{eq:model} give
\[
t:=\cT V\le C\eta_0-cB_0d_0,\qquad D^2V[e,e]\le CB_0.
\]
Here the first inequality uses $\cT(d_0^2/2)=d_0\cT d_0$;
the second follows from the uniform bounds for $d_0$ and its first two
derivatives. We only need to estimate the density where $D^2V\in\Ga$.
At such a point, $0<t\le C\eta_0$. Writing the Hessian in blocks relative
to $e^\perp\oplus\R e$, with tangential block $B$ and mixed block $z$, gives
\begin{equation}\label{eq:barrier-density}
\sigma_2(D^2V)=tD^2V[e,e]+\sigma_2(B)-|z|^2
\le C(B_0\eta_0+\eta_0^2),
\end{equation}
because $\sigma_2(B)\le\frac{n-2}{2(n-1)}(\tr B)^2$.

Choose $0<\delta_0\le1$ small enough that the right-hand side of
\eqref{eq:barrier-density} is less than $m^2$ whenever
$B_0\eta_0=\delta_0$ and $B_0\ge B_*$. Then choose $h$ small enough that
$\epsilon(h)\le2\sqrt{c\delta_0}$, where $c$ is the constant below.
For every such pair $B_0,\eta_0$, on the physical boundary we have
\[
V-\phi\ge cB_0s^4+\eta_0s^2-\epsilon(h)s^3\ge0,
\qquad
cB_0s^4+\eta_0s^2\ge2\sqrt{c\delta_0}\,s^3.
\]
On the artificial boundary $d_0\ge ch^2$. Since $\Delta v>0$, the
maximum principle bounds the normalized $v$ from above by a fixed
constant. Hence we can choose $B_0\ge B_*$ large enough that $V\ge v$
there, with $\eta_0=\delta_0/B_0$. These choices are uniform in the
boundary point and use only the stated third-derivative moduli.

If $v-V$ had a positive interior maximum, then
$D^2V\ge D^2v\in\Ga$, and ellipticity would give
$\sigma_2(D^2V)\ge\sigma_2(D^2v)\ge m^2$, a contradiction.
Thus $v\le V$ with equality at the origin. Taking the inward derivative,
$v_n(0)\le q_0-\eta_0$, and hence
\[
\tr\bigl(D^2v|_{T_0\partial\Om}\bigr)
=\Delta'\phi_2-H_0v_n(0)\ge H_0\eta_0\ge h_0\eta_0.
\]
This proves \eqref{eq:tracegap}.
\end{proof}

\subsection{Weighted quadratic comparison}
\label{wp:weighted}

Choose a collar cutoff $0\le\zeta\le1$, equal to one near $\partial\Om$,
whose support is contained in the region where $\rho=-r+Jr^2$.
Let $\Pi=I-Dr\otimes Dr$ in this collar. For a function $v$ with trace
$\varphi$, put
\[
Z_v=\zeta\frac{|\Pi D(v-\Phi)|^2}{r^2},
\]
and set $Z_v=0$ outside the collar. For smooth $v$, this quotient extends
continuously to the boundary.

We use a weighted form of the quadratic comparison identity in
\cite[Proposition~3.8]{CZZ}; see also \cite[Sections~2--3]{Lyu}.
Set $d\mu=d^{1/2}\,dx$. For smooth admissible functions
$u,v$ with the same boundary value, put
\[
\begin{gathered}
T=\Delta u,\qquad S=\Delta v,\qquad g=G(D^2u),\qquad h=G(D^2v),\\
W=u-v,\qquad \cD=TS-D^2u:D^2v-2gh.
\end{gathered}
\]
The quadratic structure gives
\begin{equation}\label{eq:defect}
\begin{split}
\cD&=\frac{TS}{2}\left|\frac{D^2u}{T}-\frac{D^2v}{S}\right|^2
+\left(g\sqrt{S/T}-h\sqrt{T/S}\right)^2\ge0,\\
\sigma_2(D^2W)&=(g-h)^2-\cD,\qquad
\cD\ge m^2(T/S+S/T)-2M^2
\end{split}
\end{equation}
when $m\le g,h\le M$.
For smooth functions $\psi,W$ vanishing on $\partial\Om$, two integrations
by parts give
\[
2\int_\Om\psi\sigma_2(D^2W)\,dx
=-\int_\Om N(D^2\psi)DW\cdot DW\,dx.
\]
Indeed $\operatorname{div}N(D^2W)=0$, and the remaining boundary terms
cancel because $DW=W_\nu\nu$.
Taking $\psi=(d+\eps)^{1/2}-\eps^{1/2}$ and then
$\eps\downarrow0$ yields
\begin{equation}\label{eq:energy}
\begin{split}
\int_\Om\cD\dm
&+\frac14\int_\Om d^{-1/2}N(D^2\rho)DW\cdot DW\,dx\\
&+\frac18\int_\Om
\frac{|D\rho|^2|DW|^2-(D\rho\cdot DW)^2}{d^2}\dm
=\int_\Om(g-h)^2\dm.
\end{split}
\end{equation}
The gradient integrals pass to the limit monotonically. They are finite
for smooth $W$ with zero boundary value.

On the support of $\zeta$, $\rho=-r+Jr^2$, so $D\rho$ is normal to the
distance level sets. Thus the last integrand in \eqref{eq:energy} is
comparable to $|\Pi DW|^2/r^2$ there. Discarding the nonnegative full-gradient term gives
\begin{equation}\label{eq:difference}
\int_\Om\cD\dm+c\int_\Om
\zeta\frac{|\Pi DW|^2}{r^2}\dm
\le\int_\Om(g-h)^2\dm.
\end{equation}

\subsection{Trace and mixed derivative estimates}
\label{wp:auxiliary}

Choose smooth $0\le\chi_\eps\le1$ with
$\chi_\eps=1$ for $r\le\eps$ and $\chi_\eps=0$ for $r\ge2\eps$,
extended by zero outside the collar, and set
\[
Z_\eps(x,p)=\zeta\frac{|\Pi(p-D\Phi)|^2}{(r+\eps)^2}.
\]
Fix a smooth nondecreasing $\eta$ with $\eta=0$ on $(-\infty,1]$ and
$\eta=1$ on $[2,\infty)$. For $b,g\in C^\infty(\overline\Om)$ with
$m\le b\le g\le M$, define
\begin{equation}\label{eq:h}
h(x,p,S)=g-(g-b)\left[\chi_\eps+(1-\chi_\eps)
\eta\!\left(\frac{S+Z_\eps(x,p)}{P}\right)\right].
\end{equation}

\begin{proposition}\label{prop:aux}
Let $\Om\subset\R^n$, $n\ge3$, be a bounded smooth domain with
$H_{\partial\Om}\ge h_0>0$, and let $\varphi\in C^\infty(\partial\Om)$.
For fixed $0<m<M$, there is a $\delta_*>0$ such that, whenever
$b,g\in C^\infty(\overline\Om)$ satisfy
$m\le b\le g\le M$ and $g-b\le\delta_*$,
the problem
\begin{equation}\label{eq:aux}
G(D^2v)=h(x,Dv,\Delta v)\quad\text{in }\Om,
\qquad v=\varphi\quad\text{on }\partial\Om
\end{equation}
with $h$ defined in \eqref{eq:h} has, for every $P\ge1$ and sufficiently
small $\eps>0$, an admissible solution $v\in C^\infty(\overline\Om)$ satisfying
\begin{equation}\label{eq:auxbound}
\Delta v+Z_v\le LP+C_b.
\end{equation}
The constant $L$ is independent of $b$ and its derivatives, $P,\eps$, and
derivatives of $g$.
\end{proposition}

\begin{proof}
We first establish the estimate. Write $S=\Delta v$. Since
$Z_\eps(x,Dv)=[r/(r+\eps)]^2Z_v$, the equation is $G(D^2v)=b$ on
each of the open sets
\begin{equation}\label{eq:frozen}
\{S>2P\},\qquad \{Z_v>8P\}.
\end{equation}
For the second assertion, use $\chi_\eps=1$ when $r\le\eps$ and
$r/(r+\eps)\ge1/2$ when $r\ge\eps$.
On these sets, use
\[
\begin{gathered}
a=DG(D^2v)=\frac{N(D^2v)}{2b},\qquad \cL=a^{ij}\partial_{ij},\\
a:D^2v=b,\qquad \tr a\ge\sqrt{\binom n2}.
\end{gathered}
\]

\emph{Gradient bound.}
Put $R_v=|Dv|^2$ and maximize $\vartheta(v)R_v$, where
$\vartheta(v)=(C_*-v)^{-1/2}$ and $C_*>\sup v+1$ is fixed using
\eqref{eq:comparison}. Boundary maxima are controlled by that estimate.
At an interior maximum with $R_v>0$,
\[
D^2v\,Dv=-\frac{R_v}{4(C_*-v)}Dv.
\]
Since $|D^2v|<S$, this gives $R_v\le4(C_*-v)S\le CS$ at the maximum.
If $S\le2P$, then $R_v\le CP$. If $S>2P$, the equation is frozen and
$a(Dv,Dv)\ge SR_v/(2M)$. At this maximum,
\[
\begin{split}
0\ge\frac{\cL(\vartheta R_v)}{\vartheta}
={}&2\sum_j a(Dv_j,Dv_j)+2Dv\cdot Db
+\frac{\vartheta'}\vartheta bR_v\\
&+\left(\frac{\vartheta''}\vartheta-
2\frac{(\vartheta')^2}{\vartheta^2}\right)R_v a(Dv,Dv).
\end{split}
\]
The last coefficient is $1/[4(C_*-v)^2]$, and $S\ge c(n,m)>0$.
Consequently $R_v^{3/2}\le C\|Db\|_\infty$ at this maximum. Since
$\vartheta$ and $\vartheta^{-1}$ are uniformly bounded, we obtain
\begin{equation}\label{eq:gradient}
\|Dv\|_\infty^2\le CP+C_b.
\end{equation}

\emph{Normalized tangential gradients.}
Fix $x_0\in\partial\Om$ and a small radius $R_0$
such that $\Om\cap B_{R_0}(x_0)$ lies where $\zeta=1$.
After a rigid motion, take $x_0=0$ and write the domain locally as
$x_n>\beta(x')$, with $\beta(0)=D\beta(0)=0$.
For each $\tau<n$, let
\[
T_\tau=\partial_\tau+
\sum_{i<n}\beta_{\tau i}(0)(x_i\partial_n-x_n\partial_i).
\]
These affine vector fields have skew-symmetric derivative matrices,
as in~\cite{ITW}.
With $\nu_{\rm in}=Dr$, their normal components satisfy
\[
(T_\tau\cdot\nu_{\rm in})(0)=
D(T_\tau\cdot\nu_{\rm in})(0)=0,
\qquad |T_\tau\cdot\nu_{\rm in}|\le C|x|^2.
\]
The tangential derivatives cancel by the choice of $T_\tau$, and the normal
derivative vanishes because $D_{\nu_{\rm in}}\nu_{\rm in}=0$.
The quadratic bound uses only the $C^3$ geometry.

Put $w=v-\Phi$ and $K=1+\|Dv\|_\infty$. On
\[
\mathcal O=\Om\cap B_{R_0}(0)
\cap\{|\Pi Dw|^2/r^2>8P\},
\]
the equation is frozen. Since $aD^2v=D^2v\,a$ and $DT_\tau$ is
skew-symmetric, differentiation gives the exact identity
\[
\cL(T_\tau w)=T_\tau b-a:D^2(T_\tau\Phi),
\qquad
|\cL(T_\tau w)|\le C(1+\|Db\|_\infty)\tr a.
\]
We compare $\pm T_\tau w$ with
\[
\mathcal B=Ad+BK|x|^2,
\qquad A=C(K+\sqrt P+1+\|Db\|_\infty).
\]
Fix $B$ sufficiently large and then choose the constant in $A$ large.
Equation~\eqref{eq:linear-barrier} gives
$\cL(\mathcal B\pm T_\tau w)\le0$ in $\mathcal O$.
On the physical boundary $Dw$ is normal, so
$|T_\tau w|\le CK|x|^2$. On the spherical boundary
$|T_\tau w|\le CK$, which is controlled by $BKR_0^2$.
On the remaining boundary, where $|\Pi Dw|^2/r^2=8P$, the normal
component estimate gives
\[
|T_\tau w|\le C\sqrt P\,r+CK|x|^2
\le Ad+BK|x|^2.
\]
Here we used $d\asymp r$ from \eqref{eq:rho}.
The minimum principle therefore yields $|T_\tau w|\le\mathcal B$ in
$\mathcal O$.

For a point $x$ in a smaller collar, choose $x_0$ to be its normal
projection. In the above coordinates $x=re_n$, and
\[
T_\tau(re_n)=e_\tau-r\sum_{i<n}\beta_{\tau i}(0)e_i
\]
is exactly tangent to the distance level set. These vectors span the
tangent plane uniformly for small fixed $R_0$. Thus the comparison and
$d\le Cr$ give $|\Pi Dw(x)|\le CAr$ whenever $x\in\mathcal O$;
outside $\mathcal O$ the defining inequality gives the same bound.
Together with \eqref{eq:gradient}, this proves
\begin{equation}\label{eq:mixed}
Z_v\le CP+C_b,\qquad
|v_{\tau\nu}|\le C\sqrt P+C_b\quad\text{on }\partial\Om.
\end{equation}
Away from this smaller collar, the first estimate follows from
\eqref{eq:gradient} because $r$ has a fixed positive lower bound.
For the boundary estimate, divide the tangential gradient bound by $r$
and let $r\downarrow0$ along normal lines; $Dw$ is normal on the boundary
and $D_{\nu_{\rm in}}\Pi=0$. Add the bounded mixed derivatives of
$\Phi$ to recover the estimate for $v$.

\emph{Trace bound.}
With $\nu$ denoting the outward normal, the identity
$D^2v|_{T\partial\Om}=D^2_{\partial\Om}\varphi+v_\nu\mathrm{II}$
and \eqref{eq:comparison} bound the tangential Hessian.
Its trace has the lower bound \eqref{eq:tracegap}. Since $h=b$ in the
boundary layer, the block identity
\[
b^2=\sigma_2(D^2v|_{T\partial\Om})
+v_{\nu\nu}\tr(D^2v|_{T\partial\Om})-|v_{T\nu}|^2
\]
and \eqref{eq:mixed} give $S\le CP+C_b$ on the boundary.
At an interior maximum of $S+|x|^2$ with $S>2P$, concavity gives
\[
\cL S\ge\Delta b,\qquad
\cL(S+|x|^2)\ge-\|\Delta b\|_\infty+\frac{(n-1)S}{b}.
\]
This bounds $S$ throughout the domain and proves \eqref{eq:auxbound}.

\emph{Solvability.}
We use the trace-dependent quotient from
\cite[Section~2.2]{Lyu}, extending the quotient calculation in
\cite[Proposition~2.2]{GZ}. For $S=\tr H>0$ set
\begin{equation}\label{eq:operator}
\cF(x,p,H)=\frac{\sigma_2(H)-h(x,p,S)^2}{S}.
\end{equation}
The function $\sigma_2(H)/S$ is concave. If $f_h=h^2$, then $h_S\le0$,
and the derivatives of the cutoff with respect to $S$ vanish unless
$S\le2P$. Direct differentiation
therefore gives
\[
S^2f_h''-2Sf_h'+2f_h
\ge2m^2-8M\|g-b\|_\infty\|\eta''\|_\infty.
\]
For fixed sufficiently small $\delta_*$, $f_h/S$ is convex and $\cF$ is
concave in $H$. Also
\[
\partial_{\lambda_i}\frac{\sigma_2(H)}S
=\frac{(S-\lambda_i)^2+\sum_{j\ne i}\lambda_j^2}{2S^2}\ge0,
\qquad D_H\cF\ge\frac{m^2}{S^2}I.
\]
The equation $\cF=0$, $S>0$, implies $H\in\Ga$.
Use the homotopy $g_t=(1-t)c+tg$, $b_t=(1-t)c+tb$, with $c\in[m,M]$.
The initial constant-density problem is classically solvable by \cite{CNS}.
Estimate \eqref{eq:auxbound} is uniform along the homotopy, so the Hessians
stay in a fixed compact subset of $\Ga$ and the linearizations are
uniformly elliptic. They have no zeroth-order term, giving openness.

To obtain a global H\"older bound for the Hessian, freeze the gradient of each
solution as a coefficient: $F_v(x,H)=\cF(x,Dv(x),H)$. The established $C^2$
bound controls the Lipschitz norms in $x$ of $F_v$ and $D_HF_v$ on compact
matrix sets. Take a convex compact $\cK\Subset\{\tr H>0\}$ containing uniform
neighborhoods of all actual Hessians and $I$, and extend by
\[
\widehat F_v(x,B)=\inf_{H\in\cK}
\{F_v(x,H)+D_HF_v(x,H):(B-H)\}.
\]
The slopes $D_HF_v(x,H)$ have common positive lower and finite upper
eigenvalue bounds, so the extension is uniformly elliptic. It is concave,
agrees with $F_v$ on $\cK$, and has spatial variation at most
$C|x-y|(1+|B|)$.
For $w=v-|x|^2/2$, normalize it as
\[
\overline F_v(x,B)=\widehat F_v(x,I+B)-F_v(x,I),
\qquad \overline F_v(x,D^2w)=-F_v(x,I).
\]
It satisfies
$|\overline F_v(x,B)-\overline F_v(y,B)|\le C|x-y||B|$:
for small $B$ integrate $D_HF_v$ on $I+tB\subset\cK$; for the remaining
$B$ absorb the constant in $1+|B|$. The right-hand side is uniformly
Lipschitz, and $w|_{\partial\Om}=\varphi-|x|^2/2$ is smooth.
Thus the concave case of \cite[Theorem~1.4(b)]{SS} gives a uniform global
$C^{2,\beta}$ bound for some $\beta>0$. Implementing the continuity method
in $C^{2,\alpha}$ with $0<\alpha<\beta$, compactness gives closedness.
Schauder bootstrapping for the original smooth equation gives smooth solutions.
These regularity constants may depend on fixed $P,\eps$ and the smooth
norms of $g$; they do not enter \eqref{eq:auxbound}.
\end{proof}

\subsection{Contraction of the trace tails}
\label{wp:tails}

\begin{proof}[Proof of Theorem~\ref{thm:wp}]
We first assume that the domain, boundary data and solution are smooth.
Write $T=\Delta u$ and $g=G(D^2u)$. Choose a smooth $b$ with
$m\le b\le g\le M$ and $\|g-b\|_\infty\le\delta\le\delta_*$,
and let $v$ solve \eqref{eq:aux}. Set $S=\Delta v$ and
\[
Y=\max(T,Z_u),\qquad E=\{S+Z_v>P\},
\qquad C_\eps=\{r<2\eps\}.
\]
By \eqref{eq:h}, $\{g\ne h\}\subset E\cup C_\eps$, and
$\mu(C_\eps)\le C\eps^{3/2}$. Consequently the right-hand side of
\eqref{eq:difference} is at most $\delta^2[\mu(E)+\mu(C_\eps)]$.

The set $E$ depends on $v$. We first control its measure
by a tail of $u$. Take a fixed sufficiently small $c_1>0$. On
$E\cap\{Y\le c_1P\}$, either $S>P/2$ or $Z_v>P/2$.
In the first case \eqref{eq:defect} gives $\cD\ge1$.
In the second, the reverse triangle inequality gives
\[
\zeta|\Pi D(u-v)|^2/r^2\ge cP.
\]
Since $P\ge1$, \eqref{eq:difference}, followed by absorption for small
fixed $\delta$, gives
\begin{equation}\label{eq:error-set}
\mu(E)\le2\mu\{Y>c_1P\}+C\delta^2\mu(C_\eps).
\end{equation}
By Proposition~\ref{prop:aux}, a threshold $P_0$ depending on
$\|b\|_{C^2}$ can be chosen so that $S,Z_v\le L_0P$ for $P\ge P_0$,
with $L_0$ independent of derivatives of $b$.
For a fixed sufficiently large $K_1$, on $\{Y>K_1P\}$ we have
\[
Y\le CP\left[\cD+
\zeta\frac{|\Pi D(u-v)|^2}{r^2}\right].
\]
If $Y=T$, use $\cD\ge m^2T/(L_0P)-2M^2$; if $Y=Z_u$,
the tangential-gradient difference is at least $Y/4$.
Integrating and using \eqref{eq:error-set} eliminate the auxiliary
solution from the estimate. Letting $\eps\downarrow0$ for each fixed
$P\ge P_0$ gives
\begin{equation}\label{eq:tail}
\int_{\{Y>K_1P\}}Y\dm
\le C_*\delta^2P\,\mu\{Y>c_1P\}.
\end{equation}
All three constants $K_1,c_1,C_*$ are independent of derivatives of $b$.

Fix $q=2p$. Multiplying \eqref{eq:tail} by $P^{q-2}$,
integrating over $P\ge P_0$, and applying Tonelli's theorem gives
\[
c_q\int_{\{Y>2K_1P_0\}}Y^q\dm
\le C_*\delta^2\int_{P_0}^{\infty}
P^{q-1}\mu\{Y>c_1P\}\,dP
\le C_q\delta^2\int_\Om Y^q\dm.
\]
Adding the integral over $\{Y\le2K_1P_0\}$ yields
\[
\int_\Om Y^q\dm
\le (2K_1P_0)^q\mu(\Om)+C_q\delta^2\int_\Om Y^q\dm.
\]
The integrals are finite at this smooth stage. Choose $\delta$ so that
$C_q\delta^2\le1/2$. We obtain
$\int_\Om Y^{2p}d^{1/2}\,dx\le C_p$. Hence
\[
\int_\Om T^p\,dx
\le\left(\int_\Om T^{2p}d^{1/2}\,dx\right)^{1/2}
\left(\int_\Om d^{-1/2}\,dx\right)^{1/2}\le C_p,
\]
since $d\asymp r$. The inequality $|D^2u|\le T$, interpolation, and the fixed amplitude bound give
the full $W^{2,p}$ estimate.

We now pass to the stated data. Set $m=\sqrt\lambda/2$ and $M=2\sqrt\Lambda$.
The order of choices is $p$, then $\delta$, then a smooth
lower approximation $b$ of $g=\sqrt f$, and finally $P_0$.
Require also $\delta\le\min(\delta_*,m/2)$.
Choose $b$ on a fixed neighborhood of $\overline\Om$ with
\[
3\delta/8\le g-b\le5\delta/8,
\qquad \|b\|_{C^2}\le C(\delta,\omega_g,\text{fixed data}).
\]
A continuous extension of $g$ with a controlled modulus and the same positive bounds, followed
by convolution and subtraction of $\delta/2$, provides such a $b\ge m$.
Thus smooth approximations of $g$ can all be compared with the same $b$.

For $C^3$ geometry and data, extend $\rho$ and $\Phi$ to a common
neighborhood and smooth them. Take the nearby zero sublevel domains
$\Om_j$ corresponding to $\Om$, and boundary values
$\varphi_j=\Phi_j|_{\partial\Om_j}$.
The boundaries converge in $C^3$; their third derivatives and those of
the boundary data have a common modulus. The mean curvature lower bound,
defining-function margin, and constants in \eqref{eq:comparison} and
\eqref{eq:tracegap} are uniform. For each $\Om_j$, reconstruct the defining
function and collar cutoff as above, with uniform constants.
Approximate $g$ on this common neighborhood
so that $b\le g_j\le b+\delta$.
The classical solutions
\[
\sigma_2(D^2u_j)=g_j^2\quad\text{in }\Om_j,\qquad
u_j=\varphi_j\quad\text{on }\partial\Om_j,
\]
exist by \cite{CNS} and satisfy the estimates just proved, uniformly in $j$.

These domains admit uniformly bounded $W^{2,p}$ extension operators.
For one $p>n$, compactness gives a subsequence converging continuously
on a neighborhood of $\overline\Om$. Its limit $u_*$ has trace $\varphi$.
Every compact subdomain of $\Om$ is eventually contained in $\Om_j$;
on these subdomains, weak continuity of Hessian measures gives
$\mu_2[u_*]=f\,dx$ by \cite[Theorem~1.1]{TW2}. Uniqueness for this Dirichlet problem identifies $u_*$
with the given solution; see \cite[Section~6]{TW2}. Weak lower semicontinuity for every finite
$p$ preserves the estimate. This proves the theorem.
\end{proof}

\section{Global \texorpdfstring{$C^{2,\alpha}$}{C2,alpha} estimates}\label{sec:c2alpha}

We fix an ambient $C^3$ extension $\Phi$ of $\varphi$ and the defining
function $d=-\rho$ from \eqref{eq:rho}. In this section, $\nu$ and $\nu_t$
denote inward unit normals; the positive mean-curvature convention remains
the outward-normal convention fixed in Section~\ref{sec:intro}.
Constants may depend on the data in Theorem~\ref{thm:c2alpha}, but not on
higher derivatives of smooth approximations. A subscript $\eps$ permits
additional dependence on $\eps$.

All differentiations below are performed on smooth admissible solutions.
Set $S=\Delta u$, $g=\sqrt f$, $a=DG(D^2u)$, and
$\cL=a^{ij}\partial_{ij}$. On $\Ga$,
\begin{equation}\label{c2:eq:ellipticity}
 a=\frac{SI-D^2u}{2g},\qquad
 \frac{g}{2S}I\le a\le\frac{S}{g}I,\qquad |D^2u|<S.
\end{equation}
Indeed, with $s_i=S-\lambda_i>0$,
$f=(S-s_i)s_i+\sigma_2(\lambda|i)\le Ss_i$.
The defining-function inequality \eqref{eq:linear-barrier} gives
$\cL d\le-\theta\tr a$ for every admissible linearization.

\subsection{Subpower bounds and smooth approximation}
\label{c2:preparation}

Theorem~\ref{thm:wp} and Sobolev embedding with $p>n/(1-\beta)$ give,
for every $0<\beta<1$,
\begin{equation}\label{c2:eq:C1beta}
 \norm{u}_{C^{1,\beta}(\overline\Om)}\le K_\beta.
\end{equation}
For smooth admissible solutions with the stated uniform data, the boundary
identity and Proposition~\ref{prop:trace} give, for the tangential Hessian
block $B$,
\begin{equation}\label{c2:eq:gap}
 |B|\le C,\qquad \tr B\ge c_0>0.
\end{equation}
\begin{lemma}\label{c2:lem:subpower}
Under the assumptions of Theorem~\ref{thm:c2alpha}, suppose that the
domain, density and boundary values are smooth. Then the admissible
solution satisfies, for every $\eps>0$,
\begin{equation}\label{c2:eq:subpower}
 |D^2u(x)|\le\Delta u(x)\le C_\eps d(x)^{-\eps}.
\end{equation}
\end{lemma}
\begin{proof}
Suppose there are $u_j,y_j$ with
$r_j^\eps\Delta u_j(y_j)\to\infty$, where $r_j$ is a fixed small
multiple of the distance to the boundary. Rescale the corresponding
interior ball to $B_2$ by
$v_j(x)=r_j^{-2}u_j(y_j+r_jx)$.
The Hessian is unchanged, and the positive density bounds and
$C^\alpha$ seminorms remain uniform.

The doubling construction and second-difference lemma
\cite[Lemmas 2.2 and 2.4, Remark 2.9]{CZZ} yield
$M_j\ge\Delta u_j(y_j)\to\infty$, points $z_j$, unit directions $e_j$,
and
\[
 h_j\asymp\frac1{M_j\log^2(e+M_j)},\qquad
 \frac{v_j(z_j+h_je_j)+v_j(z_j-h_je_j)-2v_j(z_j)}{h_j^2}
 \ge cM_j.
\]
Their constants depend on the density data, not on the gradient size
of the rescaled functions. From \eqref{c2:eq:C1beta},
\[
 M_j^\beta\le CK_\beta r_j^{\beta-1}
                       \log^{2(1-\beta)}(e+M_j).
\]
Since $r_j^\eps M_j\to\infty$, the right side is at most
$C M_j^{(1-\beta)/\eps}\log^{2(1-\beta)}(e+M_j)$ for large $j$.
Choosing $\beta>(1+\eps)^{-1}$ is a contradiction.
\end{proof}

For $0<\eps<1$, \eqref{c2:eq:ellipticity} and \eqref{c2:eq:subpower} imply
in a fixed collar that
\begin{equation}\label{c2:eq:normal-ellipticity}
 a(Dd,Dd)\ge c_\eps d^\eps.
\end{equation}
Choose $0\le\chi\le1$ smooth, equal to one on $[0,1]$, and zero
on $[2,\infty)$. For $J\ge0$, define $\psi_{J,t}:[0,\infty)\to[0,\infty)$ by
\begin{equation}\label{c2:eq:psi}
 \psi_{J,t}(0)=0,\qquad
 \psi_{J,t}'(s)=C_\eps J\int_s^{2t}r^{-\eps}\chi(r/t)\,dr
 \quad(0<s<2t),
\end{equation}
and extend it constantly after $2t$. Then
\begin{equation}\label{c2:eq:psi-bounds}
 \begin{split}
 &\cL[\psi_{J,t}(d)]
   =\psi'_{J,t}(d)\cL d+\psi''_{J,t}(d)a(Dd,Dd)
   \le-J\chi(d/t),\\
 &0\le\psi_{J,t}(d)\le C_\eps Jt^{1-\eps}\min\{d,t\},\qquad
 \psi'_{J,t}(0)\le C_\eps Jt^{1-\eps}.
 \end{split}
\end{equation}
This barrier is $C^1$ up to the physical boundary and $C^2$ in the
interior; the singular second derivative at $d=0$ causes no difficulty
for comparison or the one-sided boundary derivative.

Suppose $u,v$ have the same boundary value, are admissible, have
common positive density bounds, and satisfy \eqref{c2:eq:subpower}.
If
\[
 |G(D^2u)-G(D^2v)|
       \le\delta\,\mathbf1_{\{d<t\}}+\delta_0,
\]
then
\begin{equation}\label{c2:eq:comparison}
 |u-v|\le C_\eps\delta t^{1-\eps}\min\{d,t\}+C\delta_0d.
\end{equation}
In particular, whenever their boundary derivatives exist,
\begin{equation}\label{c2:eq:normal-comparison}
 \norm{\partial_\nu u-\partial_\nu v}_{L^\infty(\partial\Om)}
 \le C_\eps\delta t^{1-\eps}+C\delta_0.
\end{equation}
A fixed multiple of the width $t$ is harmless.
To prove these bounds, concavity at $D^2u$ gives
$\cL_u(v-u)\ge G(D^2v)-G(D^2u)$.
Use the barrier \eqref{c2:eq:psi}, and add $C\delta_0d$.
The latter handles the constant source by \eqref{eq:linear-barrier} and the
positive lower bound for $\tr a$. Thus
$\cL_u(v-u-\psi_{\delta,t}(d)-C\delta_0d)\ge0$.
The maximum principle gives the upper bound for $v-u$. Interchange
$u,v$ for the other bound. Divide by the distance and pass to the
boundary for \eqref{c2:eq:normal-comparison}.

From now on fix
\begin{equation}\label{c2:eq:epsilon-choice}
 0<\eps<\frac1{100}\min\{\alpha,1-\alpha\}.
\end{equation}
Constants may depend on this choice, and henceforth
$\beta=1-\eps/2$. Only one sufficiently large finite $p$ in
Theorem~\ref{thm:wp} is needed.

\emph{Smooth approximation.}
Extend the fixed $C^3$ functions $\rho$ and $\Phi$ to a neighborhood
of $\overline\Om$, preserving $D^2\rho-\theta I\in\Ga$
after decreasing $\theta>0$ if necessary.
Let $\omega_3$ be a common modulus for their third derivatives.

Choose an even nonnegative kernel $\eta\in C_c^\infty(B_{1/2})$ of mass
one and set
\[
 \zeta(x)=\frac{4\eta(x)-2^{-n}\eta(x/2)}3.
\]
This kernel has mass one; evenness and scaling cancel its moments of
orders one through three. It is used only for the geometry and boundary
data. Set
$\rho_t=\rho*\zeta_t$ and $\Phi_t=\Phi*\zeta_t$. Taylor's formula and
cancellation give
\begin{equation}\label{c2:sm:mollification}
\begin{split}
 \norm{D^j(\rho_t-\rho)}_\infty+
 \norm{D^j(\Phi_t-\Phi)}_\infty
 &\le Ct^{3-j}\omega_3(t),\qquad 0\le j\le2,\\
 \norm{D^4\rho_t}_\infty+\norm{D^4\Phi_t}_\infty
 &\le Ct^{-1}\omega_3(t).
\end{split}
\end{equation}
Their $C^3$ norms and third-derivative moduli are uniform in $t$.
Let $\Om_t$ be the domain bounded by the zero level of $\rho_t$
corresponding to $\Om$. For small $t$, its boundary is smooth,
$D^2\rho_t-\tfrac12\theta I\in\Ga$, and
$H_{\partial\Om_t}\ge h_0/2$. The boundary charts and the
restrictions of $\Phi_t$ have uniform $C^3$ norms and third-derivative
moduli. Their $C^4$ norms are bounded by $Ct^{-1}$.

For later comparison, choose a fixed smooth vector field transverse to
$\partial\Om$ and use its trajectories to map $\partial\Om$ onto
$\partial\Om_t$. The implicit function theorem and
\eqref{c2:sm:mollification}, followed by a cutoff in $d/t$, give a $C^3$
diffeomorphism $F_t:\overline\Om\to\overline\Om_t$ such that
\begin{equation}\label{c2:sm:diffeomorphism}
\begin{gathered}
 F_t=\mathrm{Id}\quad\text{if }d\ge t/4,\\
 \norm{F_t-\mathrm{Id}}_\infty\le Ct^3,\qquad
 \norm{DF_t-I}_\infty\le Ct^2,\\
 \norm{D^2F_t}_\infty\le Ct,\qquad
 \norm{D^3F_t}_\infty\le C.
\end{gathered}
\end{equation}
The support may be decreased by a fixed factor so that it lies where
the defining distance $d_t=-\rho_t$ is less than $t$.
The smooth transverse field avoids differentiating the nearest-point
normal field beyond its $C^2$ regularity.

Extend $g=\sqrt f$ to a fixed neighborhood, preserving the bounds
$m:=\sqrt\lambda\le g\le M:=\|f\|_{L^\infty(\Om)}^{1/2}$ and a controlled $C^\alpha$
norm. Put
$s_t=t^{4/\alpha}$ and $d_t=-\rho_t$, which is positive
in $\Om_t$ and comparable to boundary distance in a uniform collar.
With $\chi=1$ on $(-\infty,2]$ and $\chi=0$ on $[3,\infty)$,
define
\begin{equation}\label{c2:sm:density}
 g_t=\chi(d_t/t)(g*\eta_t)
       +\bigl(1-\chi(d_t/t)\bigr)(g*\eta_{s_t}),
 \qquad f_t=g_t^2.
\end{equation}
Outside the collar the cutoff is understood to be zero. These functions
are smooth and satisfy
\begin{equation}\label{c2:sm:density-bounds}
\begin{gathered}
 m\le g_t\le M,\qquad [g_t]_{C^\alpha}\le C,\qquad
 \abs{D^jg_t}\le Ct^{\alpha-j}\quad(d_t<2t, j=1,2),\\
 \abs{g_t-g}\le Ct^\alpha\quad\text{everywhere},\qquad
 \abs{g_t-g}\le Ct^4\quad(d_t\ge3t).
\end{gathered}
\end{equation}
For the uniform H\"older bound, the cutoff has Lipschitz constant
$C/t$, while the two convolutions differ by at most $Ct^\alpha$.

The classical Dirichlet theorem~\cite{CNS} gives a smooth admissible solution
\begin{equation}\label{c2:sm:classical-solution}
 G(D^2u_t)=g_t\quad\text{in }\Om_t,\qquad
 u_t=\Phi_t\quad\text{on }\partial\Om_t.
\end{equation}
Theorem~\ref{thm:wp} applies uniformly to this family.

For comparison on the fixed domain, define the corrected pullbacks as follows.
Choose a collar cutoff $\xi_t$ equal to one near the boundary and zero
when $d\ge t/4$, and set
\[
 e_t=\xi_t\bigl(\Phi_t\circ F_t-\Phi\bigr),\qquad
 U_t=u_t\circ F_t-e_t.
\]
The cutoff is fixed as a function of $d/t$. Equations
\eqref{c2:sm:mollification}--\eqref{c2:sm:diffeomorphism} imply
\begin{equation}\label{c2:sm:correction}
 U_t=\varphi\ \text{on }\partial\Om,\qquad
 \norm{D^je_t}_\infty\le Ct^{3-j}\quad(0\le j\le3).
\end{equation}

\begin{lemma}\label{c2:lem:preliminary}
The smooth approximations constructed above satisfy
\begin{equation}\label{c2:eq:preliminary-global}
 M_t:=\sup_{\Om_t}\Delta u_t\le C_\eps t^{-\eps}.
\end{equation}
\end{lemma}
\begin{proof}
Fix a boundary point, take it to be the origin, and choose coordinates
such that the inward normal is $e_n$ and the domain is locally
$x_n>b(x')$, with $b(0)=Db(0)=0$. All charts have a common radius
$R>0$ and uniform $C^3$ bounds. Suppress the subscript $t$ on the
solution, boundary extension, and linearized operator, and set
$w=u-\Phi$. For a unit tangential direction, taken to be $e_1$, put
\begin{equation}\label{c2:eq:killing-field}
 k_i=b_{1i}(0),\qquad
 T=\partial_1+\sum_{i<n}k_i(x_i\partial_n-x_n\partial_i),
 \qquad Z=Tw.
\end{equation}
The affine field $T$ generates rigid motions. Since $G$ is invariant
under these motions, differentiating the equation gives the exact identity
\begin{equation}\label{c2:eq:killing-first}
 \cL(Tu)=Tg_t.
\end{equation}
Since $T$ has affine coefficients and $\Phi_t$ has a uniform $C^3$
bound, \eqref{c2:sm:density-bounds} implies
\[
 |\cL Z|\le J+C\tr a\quad(d_t<t),\qquad J=Ct^{\alpha-1}.
\]

The field $T$ is tangent to the boundary to first order at the origin.
More precisely, since $w=0$ on the graph,
\begin{equation}\label{c2:eq:killing-defect}
 Z=\delta w_n,\qquad
 \delta(x')=k\cdot x'-b_1(x')+b(x')k\cdot Db(x'),
 \qquad |\delta(x')|\le C|x'|^2.
\end{equation}
Here the last bound uses only the uniform $C^3$ geometry. Thus
$|Z|\le C|x'|^2$ on the physical boundary. Projection to that boundary,
together with \eqref{c2:eq:C1beta}, gives
\[
 |Z(x)|\le C d_t(x)^\beta+C|x|^2
 \quad\text{in }\Om_t\cap B_R\cap\{d_t<t\}.
\]
Compare both signs of $Z$ with
\[
 P=A d_t+B|x|^2+\psi_{J,t}(d_t),\qquad
 A=C(1+t^{\beta-1}),
\]
where $B$ is a fixed sufficiently large constant. Choose $B$ first,
then the constant in $A$. On the physical boundary, $B|x|^2$ controls
$|Z|$. On $d_t=t$, the terms $A t+B|x|^2$ control the preceding bound,
and on the spherical side $B R^2$ controls the uniform bound for $Z$.
The defining-function and source-barrier inequalities give
\[
 \cL(P\pm Z)
 \le(-\theta A+2B+C)\tr a-J+J\le0.
\]
The maximum principle yields $|Z|\le P$. At the origin $Z=P=0$,
and $\partial_nZ(0)=w_{1n}(0)$ because $Dw(0)$ is normal and
$\partial_nT$ is tangential. Taking inward derivatives and using
$\psi'_{J,t}(0)\le C_\eps Jt^{1-\eps}$ proves
\begin{equation}\label{c2:eq:preliminary-mixed}
 |D^2u_t(\tau,\nu_t)|
 \le C_\eps(1+t^{\beta-1}+t^{\alpha-\eps})
 \quad\hbox{on }\partial\Om_t.
\end{equation}

Write the boundary Hessian in tangential--normal blocks as
\[
 D^2u_t=\begin{pmatrix}B&z\\z^T&D^2u_t(\nu_t,\nu_t)\end{pmatrix}.
\]
The equation gives the exact identity
\begin{equation}\label{c2:eq:normal-formula}
 D^2u_t(\nu_t,\nu_t)=\frac{|z|^2+g_t^2-\sigma_2(B)}{\tr B}.
\end{equation}
Using \eqref{c2:eq:gap} and the fixed choice $\beta=1-\eps/2$ in
\eqref{c2:eq:preliminary-mixed}, we obtain
$\sup_{\partial\Om_t}\Delta u_t\le C_\eps t^{-\eps}$.

We record the collar maximum estimate also needed below. Concavity and
\eqref{c2:sm:density-bounds} give
\begin{equation}\label{c2:eq:trace-differentiation}
 \cL_t\Delta u_t\ge\Delta g_t\ge-Ct^{\alpha-2}
 \quad(d_t<t).
\end{equation}
The maximum principle applied to
$\Delta u_t-\psi_{Ct^{\alpha-2},t}(d_t)$ in the collar yields
\begin{equation}\label{c2:eq:trace-localization}
 M_t\le
 \max\left\{\sup_{\partial\Om_t}\Delta u_t,
             \sup_{\{d_t\ge t\}}\Delta u_t\right\}
       +C_\eps t^{\alpha-\eps}.
\end{equation}
Indeed, the source barrier is nonnegative, vanishes on the physical
boundary, and has supremum at most $C_\eps t^{\alpha-\eps}$.
Lemma \ref{c2:lem:subpower} bounds the second term in the maximum by
$C_\eps t^{-\eps}$. Together with the boundary estimate, this proves
\eqref{c2:eq:preliminary-global}.
\end{proof}

\subsection{Boundary mixed derivative bounds}
\label{c2:boundary-comparison}

Set $q_t=\partial_{\nu_t}(u_t-\Phi_t)$ on $\partial\Om_t$.
We use the approximation bounds \eqref{c2:sm:mollification} and
\eqref{c2:sm:density-bounds}, together with
$1+\|D^2u_t\|_\infty\le CM_t\le C_\eps t^{-\eps}$ from
Lemma~\ref{c2:lem:preliminary}.

\begin{lemma}\label{c2:bnd:semiconcavity}
With the fixed choice \eqref{c2:eq:epsilon-choice},
\begin{equation}\label{c2:bnd:semiconcavity-bound}
 \nabla^2_{\partial\Om_t}q_t\le C_\eps t^{-1-\eps}
 g_{\partial\Om_t}.
\end{equation}
\end{lemma}

\begin{proof}
Fix a boundary point and take coordinates in which it is the origin,
$\nu_t(0)=e_n$, and the domain is locally $x_n>b(x')$, with
$b(0)=Db(0)=0$. A ball of fixed radius $R$ lies in this chart.
Suppress the subscript $t$ on $u_t,g_t,\Phi_t,d_t$ and the
linearized operator. Set $w=u-\Phi$ and
$\psi_J=\psi_{J,t}$, with $\psi_{J,t}$ as in \eqref{c2:eq:psi}.

For a unit tangential direction $e_1$, use the field $T$ and
coefficients $k_i$ from \eqref{c2:eq:killing-field}, and put $W=T^2w$.
Twice differentiating along the rigid motions generated by $T$ and using
the concavity of $G$ gives
\begin{equation}\label{c2:bnd:killing-equation}
 \cL(T^2u)
 =T^2g-G^{ij,kl}(D^2u)(Tu)_{ij}(Tu)_{kl}\ge T^2g.
\end{equation}
Since the coefficients of $T$ are affine and uniformly bounded on the
chart, \eqref{c2:sm:density-bounds} implies
\begin{equation}\label{c2:bnd:W-equation}
 \cL W\ge-J-C\norm{\Phi}_{C^4}\tr a,
 \qquad J=Ct^{\alpha-2}.
\end{equation}
On the boundary graph put
\[
 \eta(x')=e_1-b(x')k,
\]
and let $\delta$ be as in \eqref{c2:eq:killing-defect}. Then
$T=T_{\tan}+\delta\partial_n$, where $T_{\tan}$ induces
$\eta\cdot D'$ in graph coordinates. Taylor's formula gives, with
$s=\abs{x'}$,
\begin{equation}\label{c2:bnd:defect-derivative}
 \eta\cdot D'\delta
 =-\sum_{j<n}b_{11j}(0)x_j
       +O\bigl((1+\norm{b}_{C^4})s^2\bigr).
\end{equation}
Since $w=0$ on the graph and $|\delta|\le Cs^2$,
\begin{equation}\label{c2:bnd:W-boundary-identity}
 W=(\eta\cdot D'\delta)w_n
      +\delta\,\eta\cdot D'(w_n|_{\partial\Om_t})
      +\delta\,\partial_n(Tw)
       \quad\hbox{on }\partial\Om_t.
\end{equation}
The last two terms are bounded in absolute value by $CM_ts^2$, and
\[
 \abs{w_n(x',b(x'))-w_n(0)}\le CM_ts.
\]
The uniform gradient bound, \eqref{c2:bnd:defect-derivative}, and
\eqref{c2:bnd:W-boundary-identity} therefore yield
\begin{equation}\label{c2:bnd:linear-boundary-error}
 \abs{W-\ell}\le C(M_t+1+\norm{b}_{C^4})s^2,
 \qquad
 \ell(x')=-w_n(0)\sum_{j<n}b_{11j}(0)x_j.
\end{equation}
Here $\cL\ell=0$ and $\partial_n\ell=0$.

In $\Om_t\cap B_R\cap\{d<t\}$ compare $W-\ell$ with
\[
 P=Ad+B\abs{x}^2+\psi_J(d),\qquad
 B=C(M_t+t^{-1}),\qquad
 A=C(M_t/t+t^{-1}).
\]
The constants are chosen in the order $B$, then $A$.
On the physical boundary, \eqref{c2:bnd:linear-boundary-error} gives
$W-\ell\le P$. On $d=t$, the estimate $\abs{W-\ell}\le CM_t$
is absorbed by $At$. On the spherical boundary, it is absorbed by
$BR^2$. By \eqref{eq:linear-barrier}, \eqref{c2:eq:psi-bounds}, and
\eqref{c2:bnd:W-equation},
\[
 \cL(W-\ell-P)
 \ge\bigl(A\theta-2B-C\norm{\Phi}_{C^4}\bigr)\tr a\ge0.
\]
The maximum principle gives $W\le\ell+P$. Equality holds at the
origin, so taking the inward derivative gives
\begin{equation}\label{c2:bnd:W-normal}
 W_n(0)\le C_\eps\bigl(M_t/t+t^{-1}+Jt^{1-\eps}\bigr)
          \le C_\eps t^{-1-\eps}.
\end{equation}

Differentiating $T^2w$ at the origin and using
$w_{i1}=-b_{i1}w_n$ there gives
\[
 W_n(0)=w_{11n}(0)
             +b_{11}(0)w_{nn}(0)
             +\abs{k}^2w_n(0).
\]
On the boundary,
\[
 \partial_{\nu_t}w
   =\sqrt{1+\abs{Db}^2}\,w_n(x',b(x')).
\]
Its second derivative in $x_1$ at the origin is the same expression.
The Christoffel symbols of the graph metric vanish at the origin.
Consequently,
\begin{equation}\label{c2:bnd:exact-normal-identity}
 W_n(0)=\nabla^2_{\partial\Om_t}
           (\partial_{\nu_t}w)(e_1,e_1).
\end{equation}
Equations \eqref{c2:bnd:W-normal} and \eqref{c2:bnd:exact-normal-identity}
prove \eqref{c2:bnd:semiconcavity-bound}.
\end{proof}

\begin{lemma}\label{c2:lem:boundary}
There exist $t_0>0$ and $C>0$, depending only on the data in
Theorem~\ref{thm:c2alpha}, such that, for $t_j=2^{-j}t_0$,
\begin{equation}\label{c2:bnd:uniform-boundary-hessian}
 \norm{D^2u_{t_j}}_{L^\infty(\partial\Om_{t_j})}\le C
 \qquad(j\ge0).
\end{equation}
\end{lemma}

\begin{proof}
\emph{Comparison on the fixed domain.}
The Hessian chain rule for the pullbacks $U_t$ gives
\[
 D^2U_t=(DF_t)^T(D^2u_t\circ F_t)DF_t
       +(Du_t\circ F_t)D^2F_t-D^2e_t.
\]
The global gradient bound, \eqref{c2:sm:diffeomorphism}, and
Lemma~\ref{c2:lem:preliminary} give
\begin{equation}\label{c2:sm:hessian-error}
 \abs{D^2U_t-D^2u_t\circ F_t}
 \le C_\eps(t+t^{2-\eps}).
\end{equation}
Since $\sigma_2$ is quadratic, its change is bounded by
$C_\eps t^{1-\eps}$. For small $t$, both the trace and $\sigma_2$ remain
positive, so $D^2U_t\in\Ga$. Consequently
\begin{equation}\label{c2:sm:transformed-density}
 \abs{G(D^2U_t)-g_t\circ F_t}
 \le C_\eps t^{1-\eps}\,\mathbf1_{\{d<t/4\}}.
\end{equation}
The support of $F_t-\mathrm{Id}$ lies inside the region where $g_t$ is
coarsely mollified. Thus its composition error is at most
$Ct^3t^{\alpha-1}=Ct^{2+\alpha}$ there.

Write $\widetilde g_t=G(D^2U_t)$. Since $1-\eps>\alpha$,
\eqref{c2:sm:density-bounds} and \eqref{c2:sm:transformed-density} imply
\begin{equation}\label{c2:sm:dyadic-density}
 \abs{\widetilde g_t-\widetilde g_{2t}}
 \le C_\eps t^\alpha\mathbf1_{\{d<Ct\}}+Ct^4.
\end{equation}
The global preliminary bound and \eqref{c2:sm:hessian-error} give,
on the larger collar $d<3Ct$,
$a^{ij}d_i d_j\ge c_\eps t^\eps$, hence
$a^{ij}d_i d_j\ge c_\eps d^\eps$ after adjusting the constant, for the
linearizations at both $U_t$ and $U_{2t}$. This collar contains the
full support of the second derivative of the source barrier.
Moreover $\cL d\le-\theta\tr a$. Apply the thin-layer comparison
estimate to \eqref{c2:sm:dyadic-density}, and add the fixed barrier
$Ct^4d$ for the error outside the collar. If
$Q_t=\partial_{\nu}(U_t-\Phi)|_{\partial\Om}$, this gives
\begin{equation}\label{c2:sm:dyadic-comparison}
 \norm{U_t-U_{2t}}_\infty\le C_\eps t^{2+\alpha-\eps},\qquad
 \norm{Q_t-Q_{2t}}_\infty\le C_\eps t^{1+\alpha-\eps}.
\end{equation}
Since $u_t-\Phi_t$ vanishes on $\partial\Om_t$, its tangential gradient
is zero. Near $\partial\Om$ the correction in \eqref{c2:sm:correction}
is exact, so
\begin{equation}\label{c2:sm:normal-transformation}
 Q_t=a_t(q_t\circ F_t),\qquad
 a_t=(DF_t\nu)\cdot(\nu_t\circ F_t).
\end{equation}
On the fixed boundary,
\[
 a_t=1+O(t^2),\qquad
 \abs{\nabla a_t}\le Ct,\qquad
 \abs{\nabla^2a_t}\le C,\qquad
 \abs{q_t}+\abs{\nabla q_t}\le C_\eps t^{-\eps}.
\]
Derivatives are taken with respect to the corresponding induced boundary
metrics. The last bound follows from Lemma~\ref{c2:lem:preliminary}
and the uniform $C^3$ data. Applying the boundary chain rule to
\eqref{c2:sm:normal-transformation}, using $a_t>0$ and
\eqref{c2:bnd:semiconcavity-bound}, gives
\begin{equation}\label{c2:sm:fixed-semiconcavity}
 \nabla^2_{\partial\Om}Q_t\le
 C_\eps t^{-1-\eps}g_{\partial\Om}.
\end{equation}
Indeed, the leading term is $a_t$ times the pullback of the intrinsic
Hessian of $q_t$. All remaining terms contain only $q_t$, its first
derivatives, and derivatives of $F_t$ of order at most three.

\emph{The boundary gradient recurrence.}
On $\Sigma=\partial\Om$, put $t_j=2^{-j}t_0$ and
$K_j=\sup_\Sigma|\nabla_\Sigma Q_{t_j}|$.
For $j\ge1$, \eqref{c2:sm:dyadic-comparison} and
\eqref{c2:sm:fixed-semiconcavity} give
\[
 \|Q_{t_j}-Q_{t_{j-1}}\|_\infty\le E_j:=C_\eps t_j^{1+\alpha-\eps},
 \qquad
 \nabla^2_\Sigma Q_{t_j}\le L_jg_\Sigma,
 \quad L_j:=C_\eps t_j^{-1-\eps}.
\]
At $p\in\Sigma$ with $\nabla_\Sigma Q_{t_j}(p)\ne0$, take the unit-speed
geodesic $\gamma$ starting in the direction of its negative gradient.
Along this curve, $(Q_{t_j}\circ\gamma)''\le L_j$. One-dimensional
Taylor's formula and the bound for $Q_{t_j}-Q_{t_{j-1}}$ give, for small
$h>0$,
\[
 -K_{j-1}h-2E_j
 \le Q_{t_j}(\gamma(h))-Q_{t_j}(p)
 \le-h|\nabla_\Sigma Q_{t_j}(p)|+\frac{L_j}{2}h^2.
\]
Taking $h=2\sqrt{E_j/L_j}\asymp t_j^{1+\alpha/2}$, which is admissible
for a fixed sufficiently small $t_0$, yields
\begin{equation}\label{c2:bnd:lipschitz-recurrence}
 K_j\le K_{j-1}+2\sqrt{E_jL_j}
       \le K_{j-1}+C_\eps t_j^{\alpha/2-\eps}.
\end{equation}
The preliminary Hessian and pullback bounds give
$K_0\le C_\eps t_0^{-\eps}$. Summing the increments gives $K_j\le C$,
with $C$ depending only on the stated data.

In \eqref{c2:sm:normal-transformation}, $a_t\ge c>0$, the $C^1$ norms
of $a_t$ are bounded, and $F_t$ and its inverse have uniformly bounded
first derivatives. The gradient bound and the estimate for $K_j$
therefore imply
\[
 \norm{\nabla_{\partial\Om_{t_j}}q_{t_j}}_{L^\infty}\le C.
\]
For a unit tangential vector $\tau$, the zero boundary value gives
\[
 \tau(q_{t_j})
 =D^2(u_{t_j}-\Phi_{t_j})(\tau,\nu_{t_j}),
\]
since $D(u_{t_j}-\Phi_{t_j})$ is normal and $D_\tau\nu_{t_j}$ is
tangential. Adding the bounded derivatives of $\Phi_{t_j}$ yields
\begin{equation}\label{c2:bnd:uniform-mixed}
 \abs{D^2u_{t_j}(\tau,\nu_{t_j})}\le C
       \quad\hbox{on }\partial\Om_{t_j}.
\end{equation}
The boundary block formula \eqref{c2:eq:normal-formula}, the tangential
bounds \eqref{c2:eq:gap}, and \eqref{c2:bnd:uniform-mixed} give
\eqref{c2:bnd:uniform-boundary-hessian}.
\end{proof}

\subsection{Global trace bounds}
\label{c2:trace-comparison}

Fix $t_0$ as in Lemma~\ref{c2:lem:boundary} and write
$t=t_j=2^{-j}t_0$, $j\ge1$.
The boundary estimate \eqref{c2:bnd:uniform-boundary-hessian}
and the trace localization \eqref{c2:eq:trace-localization} give
\begin{equation}\label{c2:eq:move-maximum}
 M_t\le
 \max\left\{C_0,\sup_{\{d_t\ge t\}}\Delta u_t\right\}
       +C_\eps t^{\alpha-\eps},
\end{equation}
with $C_0$ depending only on the stated data.
Either $M_t\le C_0+C_\eps t^{\alpha-\eps}$, or there is $x_t$ with
$d_t(x_t)\ge t$ and $\Delta u_t(x_t)\ge M_t-C_\eps t^{\alpha-\eps}$.
Both $F_t,F_{2t}$ equal the identity and $e_t,e_{2t}$ vanish where
$d\ge t/2$. Since $|d_t-d|\le Ct^3$, a sufficiently small fixed $c_1>0$
puts $B_R(x_t)$, $R=c_1t$, inside this common region.

The fixed-domain comparison therefore gives
\begin{equation}\label{c2:eq:interior-difference}
 \norm{u_t-u_{2t}}_{L^\infty(B_R(x_t))}
       \le C_\eps t^{2+\alpha-\eps}.
\end{equation}

\begin{lemma}\label{c2:lem:finite-trace}
The approximations constructed above satisfy
\[
 M_t\le\max\{C_0,M_{2t}\}+Ct^{\alpha/10}.
\]
\end{lemma}

\begin{proof}
By \eqref{c2:eq:move-maximum}, only the alternative with $x_t$
requires proof. Put $M=M_t$; the lower density bound gives
$M\ge c(n,\lambda)>0$. Translate $x_t$ to the origin. Choose
$0\le\chi_1\le1$ equal to one on $B_{3R/4}$ and supported in $B_R$, and
$0\le\chi_0\le1$ supported in $B_{R/4}$, with $\chi_0(0)=1$ and
\[
 \frac{|D\chi_0|^2}{\chi_0}+|D^2\chi_0|\le CR^{-2}
 \quad\hbox{where }\chi_0>0.
\]
Subtract the affine tangent $\ell$ of $u_t$ at the origin, and set
\[
 \begin{gathered}
 \widetilde u=\chi_1(u_t-\ell),\qquad R_s=(I-s\Delta)^{-1},\\
 H_s=R_sD^2\widetilde u,\qquad T_s=\tr H_s,\qquad
 N(s)=\max_x\chi_0(x)T_s(x).
 \end{gathered}
\]
Here $R_s$ acts on $\R^n$, after extension by zero. Since
$|D^2u_t|\le\Delta u_t\le M$ on $B_R$, we have
$|D^2\widetilde u|\le CM$ and
\[
 N(0)\ge M-e,\qquad e=C_\eps t^{\alpha-\eps}.
\]
For small $t$, $e\le M/8$.

Write $r=\sqrt{s}$. The convolution kernel $G_s$ of $R_s$ has
the representation \cite[Lemma~2.5]{CZZ}
\[
 G_s(z)=\int_0^\infty e^{-\tau}K_{s\tau}(z)\,d\tau,
 \qquad K_a(z)=(4\pi a)^{-n/2}e^{-|z|^2/(4a)},\quad a>0.
\]
It is positive and has total mass one. The same lemma gives,
uniformly for $x\in B_{R/4}$,
\[
 \int_{\R^n\setminus B_{3R/4}}G_s(x-y)\,dy
 \le\epsilon_r:=C(1+R/r)^N e^{-cR/r}.
\]
We prove that, whenever $N(s)\ge M/2$ and
\[
 A(r):=r^\alpha+\frac{M^2r^2}{R^2}+M^2\epsilon_r
\]
is sufficiently small,
\begin{equation}\label{c2:eq:resolvent-log-derivative}
 \left|s\frac{d}{ds}\log N(s)\right|\le C\sqrt{A(r)}
 \quad\hbox{for almost every }s.
\end{equation}
\emph{Resolvent estimate.}
Fix a maximum point $x_s$ of $\chi_0T_s$. Then
$T_s(x_s)\asymp M$, $\chi_0(x_s)\ge c$, and the first and second
derivative tests give
\[
 D^2T_s(x_s)\le CMR^{-2}I.
\]
All averages below are centered at $x_s$. Put $I=B_{3R/4}$ and define
\[
 \begin{gathered}
 b=\int_I G_s(x_s-y)\,dy,\qquad
 d\bar\nu=b^{-1}\mathbf1_I(y)G_s(x_s-y)\,dy,\\
 \bar H=\int_I D^2u_t\,d\bar\nu,\qquad \bar T=\tr\bar H.
 \end{gathered}
\]
Thus $1-b\le\epsilon_r$ and
\[
 |H_s(x_s)-\bar H|\le CM\epsilon_r.
\]
Concavity of $G=\sqrt{\sigma_2}$ gives $\bar H\in\Ga$ and
$\sigma_2(\bar H)\ge\lambda$. Hence $H_s(x_s)\in\Ga$ when
$M^2\epsilon_r$ is small. The resolvent identity yields
\[
 H_s(x_s)-D^2u_t(x_s)=sD^2T_s(x_s)\le CMr^2R^{-2}I.
\]
The Newton tensors along the segment joining these admissible matrices
are positive and have trace at most $CM$. Thus
\[
 \sigma_2(H_s(x_s))\le f_t(x_s)+CM^2r^2R^{-2}.
\]
Let $S=\Delta u_t$, $g=\sqrt{f_t}$, and $m=\int_I g\,d\bar\nu$.
Concavity of $\sigma_2/\sigma_1$ and Cauchy--Schwarz give
\[
 m^2\le\bar T\int_I\frac{f_t}{S}\,d\bar\nu
       \le\sigma_2(\bar H).
\]
Using $|m^2-f_t(x_s)|\le Cr^\alpha$ and the preceding matrix estimates,
the quotient-concavity argument in \cite[equation (2.14)]{CZZ} gives
\begin{equation}\label{c2:eq:resolvent-defect}
 0\le\bar T\int_I\frac{f_t}{S}\,d\bar\nu-m^2\le CA(r).
\end{equation}

Expanding the square in \eqref{c2:eq:resolvent-defect} gives
\[
 \bar T\int_I\frac{(g-mS/\bar T)^2}{S}\,d\bar\nu
 =\bar T\int_I\frac{f_t}{S}\,d\bar\nu-m^2\le CA(r).
\]
Since $S\le M$, $\bar T\asymp M$, and $m$ is bounded below by a
positive constant, it follows that
\[
 \int_I\left|\frac{S}{\bar T}-\frac gm\right|^2d\bar\nu
 \le CA(r).
\]
Moreover,
$\int_I|g/m-1|^2\,d\bar\nu\le Cr^{2\alpha}$ by the kernel moment
bound and the $C^\alpha$ norm of $g$. Together with
$|T_s(x_s)-\bar T|\le CM\epsilon_r$, these estimates imply
\[
 \int_I\left|\frac{S}{T_s(x_s)}-1\right|^2d\bar\nu
 \le C\bigl(A(r)+r^{2\alpha}+\epsilon_r^2\bigr).
\]
On $\R^n\setminus I$ the corresponding quotient with
$\Delta\widetilde u$ is bounded. Since $M$ is bounded below and $r<1$, the full
kernel average consequently satisfies
\begin{equation}\label{c2:eq:resolvent-quadratic-concentration}
 R_s\!\left[\left(\frac{\Delta\widetilde u}{T_s(x_s)}-1\right)^2
       \right](x_s)\le CA(r).
\end{equation}

For $w\in C_c^\infty(\R^n)$, the kernel representation and
$s\partial_sK_{s\tau}=\tau\partial_\tau K_{s\tau}$ give, by
integration by parts in $\tau$,
\[
 s\partial_sR_sw
 =\int_0^\infty(\tau-1)e^{-\tau}(K_{s\tau}*w)\,d\tau.
\]
The boundary terms vanish because $K_{s\tau}*w$ is bounded.
Since $\int_{\R^n}K_a=1$,
$\int_0^\infty(\tau-1)e^{-\tau}\,d\tau=0$, and
$\int_0^\infty(\tau-1)^2e^{-\tau}\,d\tau=1$,
Cauchy--Schwarz gives, for every constant $c$,
\begin{equation}\label{c2:eq:resolvent-square}
 |s\partial_sR_sw(x)|^2\le R_s[(w-c)^2](x).
\end{equation}
Apply this to the fixed function $w=\Delta\widetilde u$, with
$c=T_s(x_s)$. Equation~\eqref{c2:eq:resolvent-quadratic-concentration}
then gives
\[
 \left|\frac{s\partial_sT_s(x_s)}{T_s(x_s)}\right|
 \le C\sqrt{A(r)}.
\]
The function $N$ is locally Lipschitz for $s>0$; wherever it is
differentiable, its derivative equals
$\chi_0(x_s)\partial_sT_s(x_s)$. This proves \eqref{c2:eq:resolvent-log-derivative}.

\emph{Comparison of successive approximations.}
Choose $\delta=\alpha/4$ and $h=t^{1+\delta}$. Integrating with
$ds/s=2\,dr/r$ gives a relative loss bounded by
\[
 I_h:=C\int_0^h\sqrt{A(r)}\,\frac{dr}{r}
 \le C\left[h^{\alpha/2}+\frac{Mh}{R}
       +M(1+R/h)^Ne^{-cR/h}\right].
\]
The bound $M\le C_\eps t^{-\eps}$ and
$\delta>\eps$ show that $A(r)$ is uniformly small for $0<r\le h$
and $I_h=o(1)$. If $N$ first reached $M/2$ before $h^2$, integration
of \eqref{c2:eq:resolvent-log-derivative} up to that first time, with
$N(0)\ge7M/8$, would give a contradiction. Thus $N(s)>M/2$ throughout
$[0,h^2]$, and
\[
 N(h^2)\ge(M-e)e^{-I_h}\ge M-e-CMI_h.
\]
In particular, at some $z\in B_{R/4}$,
\[
 R_{h^2}\Delta\widetilde u(z)\ge M-e-CMI_h.
\]

Use the same affine function and cutoff to set
$\widetilde v=\chi_1(u_{2t}-\ell)$ and
$\widetilde w=\widetilde u-\widetilde v=\chi_1(u_t-u_{2t})$.
Then
\[
 |R_{h^2}\Delta\widetilde w|
 =h^{-2}|R_{h^2}\widetilde w-\widetilde w|
 \le2h^{-2}\|u_t-u_{2t}\|_\infty.
\]
Positivity of the kernel and $\Delta u_{2t}\le M_{2t}$ on $B_{3R/4}$
give, at the same point $z$,
\[
 R_{h^2}\Delta\widetilde v(z)
 \le M_{2t}+C(M_{2t}+R^{-2})(1+R/h)^Ne^{-cR/h}.
\]
The stated $C^1$ bounds suffice for the cutoff error here. Combining
these inequalities gives
\begin{align*}
 M_t\le M_{2t}+C_\eps\bigl[
 &t^{\alpha-\eps}
 +t^{(1+\delta)\alpha/2-\eps}
 +t^{\delta-2\eps}
 +t^{\alpha-\eps-2\delta}\bigr]+O(t^N)
\end{align*}
for every fixed $N$. All four powers exceed $\alpha/10$ when
$\delta=\alpha/4$ and $0<\eps<\alpha/40$, as ensured by
\eqref{c2:eq:epsilon-choice}. The bounded alternative in
\eqref{c2:eq:move-maximum} obeys the same conclusion.
\end{proof}

\subsection{The \texorpdfstring{$C^{2,\alpha}$}{C2,alpha} estimate}
\label{c2:holder}

\begin{proof}[Proof of Theorem~\ref{thm:c2alpha}]
Lemma~\ref{c2:lem:finite-trace} and the initial bound
$M_{t_0}\le C_\eps t_0^{-\eps}$ give, by summation,
\begin{equation}\label{c2:eq:uniform-global-trace}
 \sup_j M_{2^{-j}t_0}\le C.
\end{equation}
The density lower bound and \eqref{c2:eq:uniform-global-trace} keep
the Hessians in a fixed compact subset of $\Ga$, where the
linearizations are uniformly elliptic by \eqref{c2:eq:ellipticity}. Set
\[
 \cK=\{H\in\Ga:G(H)\ge\sqrt{\lambda}/2,\ \tr H\le2C\},
 \qquad
 \widehat G(B)=\inf_{H\in\cK}DG(H):B.
\]
Concavity and homogeneity give $\widehat G=G$ on $\cK$. The common
positive lower and finite upper eigenvalue bounds for $DG(H)$ make
$\widehat G$ globally uniformly elliptic. It is concave, vanishes at
zero, and agrees with $G$ on a fixed neighborhood of all the attained
Hessians.
The global $C^{2,\beta_0}$ theorem gives a uniform bound for some
$0<\beta_0<\alpha$; apply \cite[Theorem 1.4(b)]{SS} to this extension.
The required domain and boundary norms are uniform, since the approximations
have common $C^3$ bounds.

We obtain the exponent $\alpha$ by linear Schauder estimates.
Choose a small number $\eta>0$. The uniform $C^{2,\beta_0}$ estimate
provides a radius $r>0$, independent of $t$, on whose interior and
boundary patches the Hessian oscillation is at most $\eta$.
Fix an attained Hessian $H_0$ in a patch and put $a_0=DG(H_0)$ and
\[
 R(H)=G(H_0)+a_0:(H-H_0)-G(H).
\]
Homogeneity gives $a_0:H_0=G(H_0)$, so in the original Euclidean
coordinates the equation is exactly
\[
 a_0:D^2u_t=g_t+R(D^2u_t).
\]
Smoothness of $G$ on the fixed Hessian neighborhood implies
\[
 [R(D^2u_t)]_{C^\alpha(\text{patch})}
       \le C\eta[D^2u_t]_{C^\alpha(\text{patch})}.
\]
Apply the scaled local interior or boundary linear Schauder estimate
for the constant matrix $a_0$; see \cite[Chapter 6]{GT}. Taking the
supremum over patches of this fixed radius, and using the Hessian
bound for pairs farther apart than a fixed multiple of $r$, gives
\[
 [D^2u_t]_{C^\alpha(\overline\Om_t)}
 \le C_r+C\eta[D^2u_t]_{C^\alpha(\overline\Om_t)}.
\]
In this scaled estimate the coefficient $C$ is independent of $r$
and $t$; the finite lower-order terms are included in $C_r$.
Choose $\eta$ so that $C\eta<1/2$. The solutions are smooth, so the
seminorm being absorbed is finite before the estimate is taken.
Thus
\[
 \norm{u_t}_{C^{2,\alpha}(\overline\Om_t)}\le C
 \quad\hbox{along }t=2^{-j}t_0.
\]

The pullbacks $U_t=u_t\circ F_t-e_t$ obey the same bound on the
fixed domain, using the uniform $C^3$ bounds for $F_t$ and $e_t$.
A subsequence converges in $C^2(\overline\Om)$, and lower
semicontinuity of the H\"older seminorm preserves the
$C^{2,\alpha}$ bound. Its density is $f$: the density-smoothing
error tends uniformly to zero, as does the pullback density error
$O(t^{1-\eps})$. Its boundary value is $\varphi$.
Admissibility and comparison identify the limit with the given
continuous $2$-convex solution.
\end{proof}

\section{Counterexamples}\label{sec:examples}

The following examples show that the global $W^{2,p}$ conclusion for every
finite $p$, and hence the global $C^{2,\alpha}$ conclusion, fail at
$C^{2,1}$ regularity, even when the density is identically one.
Analogous failures for the Monge--Amp\`ere equation were established by Wang~\cite{Wang}.
Both examples are obtained by extending a two-dimensional
Monge--Amp\`ere solution constantly in the remaining variables.
In both examples $n\ge3$, $x'=(x_1,\ldots,x_{n-1})$, and
\[
 t=x_n+x_1^2,\qquad U(x)=\frac23t^{3/2}\quad(x_n>0).
\]
The only nonzero block of its Hessian, in the $(x_1,x_n)$ coordinates, is
\[
 \begin{pmatrix}
 2\sqrt t+2x_1^2/\sqrt t&x_1/\sqrt t\\
 x_1/\sqrt t&1/(2\sqrt t)
 \end{pmatrix}.
\]
Its determinant equals one and both diagonal entries are positive. Thus
$D^2U\in\Ga$ and $\sigma_2(D^2U)=1$. On bounded sets with $t>0$ small,
\begin{equation}\label{eq:ce-hessian}
 c\,t^{-1/2}\le |D^2U|\le C\,t^{-1/2}.
\end{equation}

\subsection{\texorpdfstring{$C^{2,1}$}{C2,1} boundary values on a ball}
Let $\Om=B_1(e_n)$ and $\varphi=U|_{\partial\Om}$. The only boundary
point where $t=0$ is the origin. Near this point the boundary is
\[
 x_n=\beta_0(x')=1-\sqrt{1-|x'|^2}
       =\frac12|x'|^2+O(|x'|^4).
\]
Consequently, writing $s=|x'|$ and
$Q(x')=x_1^2+\beta_0(x')$, we have $Q\asymp s^2$, $DQ=O(s)$,
and bounded derivatives of $Q$ through order three. Differentiation of
$\varphi(x',\beta_0(x'))=\frac23Q^{3/2}$ gives, for $s>0$,
\[
 |D^j(\varphi(x',\beta_0(x')))|\le Cs^{3-j},\qquad 0\le j\le3.
\]
The second derivatives extend continuously by zero at the origin and are
Lipschitz. Hence $\varphi\in C^{2,1}(\partial\Om)$. Along the $x_1$ axis,
\[
 \varphi(ae_1,\beta_0(ae_1))
 =\frac23\Bigl(\frac32\Bigr)^{3/2}|a|^3+O(|a|^5),
\]
so this trace is not $C^3$.

Inside the ball, $x_n\le t\le3x_n$. At height $x_n=z$, the horizontal
section has $(n-1)$-dimensional measure
\[
 \omega_{n-1}(2z-z^2)^{(n-1)/2}\asymp z^{(n-1)/2}
 \qquad (0<z\ll1).
\]
Since $U$ and $DU$ are bounded, \eqref{eq:ce-hessian} yields the exact threshold
\begin{equation}\label{eq:ce-threshold}
 U\in W^{2,p}(\Om)
 \quad\Longleftrightarrow\quad p<n+1,
 \qquad 1\le p<\infty.
\end{equation}
At $p=n+1$ the Hessian integral diverges logarithmically.

This example also explains the dependence on the third-derivative modulus
of the boundary data. The trace $\varphi$ admits smooth approximations
with uniformly bounded $C^3$ norms, and comparison implies that the
corresponding solutions with density one converge uniformly to $U$.
A $W^{2,p}$ estimate depending only on those norms would therefore give
$U\in W^{2,p}(\Om)$ by weak compactness, contradicting
\eqref{eq:ce-threshold} whenever $p\ge n+1$.

\subsection{A \texorpdfstring{$C^{2,1}$}{C2,1} domain with polynomial data}
Set
\[
 \Phi(x)=x_n-\frac12|x|^2,\qquad
 \Om_* =\{x_n>0:U(x)<\Phi(x)\}.
\]
Put
$F=U-\Phi$. In the upper half-space,
\[
 D^2F=I+D^2U\ge I.
\]
Thus $\Om_*$ is convex. It is nonempty, since $F(ze_n)<0$ for small
$z>0$, and it is bounded because $U<\Phi$ implies
$\Om_*\subset B_1(e_n)$. Its closure meets $\{x_n=0\}$ only at the origin.
At every other boundary point $DF\ne0$: a critical point of the convex
function $F$ would be a global minimum, whereas $F$ takes negative values.
The boundary is therefore smooth away from the origin.

Near the origin its lower boundary is a graph $x_n=\beta(x')$ satisfying
\begin{equation}\label{eq:ce-graph}
 \beta=\frac12|x'|^2+\frac12\beta^2
               +\frac23(x_1^2+\beta)^{3/2}.
\end{equation}
For $s=|x'|>0$ sufficiently small,
\[
F(x',s^2/2)>0>F(x',s^2),\qquad
F_{x_n}=\sqrt{x_1^2+x_n}-1+x_n<-1/2
\]
on this interval. Thus the equation has a unique root between $s^2/2$
and $s^2$. Set $\beta(0)=0$ and write
$x'=s\theta$ and $\beta=s^2B(s,\theta)$. Then
\[
 B=\frac12+\frac23s(\theta_1^2+B)^{3/2}+\frac12s^2B^2.
\]
At $s=0$ the derivative of this equation with respect to $B$, with all
terms on the left, is one. The implicit function theorem gives a smooth
$B$ for small $s\ge0$, uniformly for $\theta\in\mathbb S^{n-2}$, and
$B(0,\theta)=1/2$. In particular,
\[
 \bigl|D^j\bigl(\beta(x')-\tfrac12|x'|^2\bigr)\bigr|
 \le Cs^{3-j},\qquad 0\le j\le3,\quad s>0.
\]
It follows that $D^2\beta\to I$ at the origin and that $\beta\in C^{2,1}$.
Moreover, \eqref{eq:ce-graph} gives
\[
 \beta(ae_1)=\frac12a^2
       +\frac23\Bigl(\frac32\Bigr)^{3/2}|a|^3+O(|a|^4),
\]
which excludes $C^3$ regularity.

The domain is uniformly convex. Indeed, $|DF|$ is bounded on its closure,
and at every smooth boundary point its second fundamental form, with the
outward normal, satisfies
\[
 \mathrm{II}(\tau,\tau)
   =\frac{D^2F[\tau,\tau]}{|DF|}\ge c|\tau|^2.
\]
At the origin the same conclusion follows from $D^2\beta(0)=I$. Hence
$\partial\Om_*\in C^{2,1}$ and its principal curvatures have a uniform
positive lower bound.

On $\partial\Om_*$ we have $U=\Phi$, including at the origin. Thus the
boundary values are the restriction of an ambient quadratic polynomial.
Finally, $\beta(x')\asymp |x'|^2$ implies that the
horizontal sections of $\Om_*$ have measure comparable to
$z^{(n-1)/2}$ as $z\downarrow0$. Also $z\le t\le3z$, since
$\Om_*\subset B_1(e_n)$. Therefore \eqref{eq:ce-threshold} holds with
$\Om_*$ in place of $\Om$.

\section*{Acknowledgements}

The author is deeply grateful to Professors Dongsheng Li and Yu Yuan for their guidance and support.
AI tools were used as an aid in preliminary mathematical exploration and in refining the exposition.
All mathematical arguments were independently verified by the author, who takes full responsibility for the manuscript.

\begingroup
\renewcommand{\baselinestretch}{1}\selectfont

\endgroup
\end{document}